# Asymptotic Solution of the Anti-Plane Problem for a Two-Dimensional Lattice

## N.I. Aleksandrova

*N.A. Chinakal Institute of Mining, Siberian Branch, Russian Academy of Sciences, Krasnyi pr. 91, Novosibirsk, 630091 Russia, e-mail: alex@math.nsc.ru*

**Abstract.** Propagation of unsteady waves under the effect of a step point load on a square lattice of spring-connected masses is investigated. The problem is solved by two methods. Asymptotic solutions at large time intervals, which describe the behavior of long-wave perturbations, are derived analytically. The solution over the whole time interval for the waves of the entire spectral range is derived by the finite difference method. These solutions are compared, and their good agreement is shown.

Previously, the theory of deforming the rock massif as a homogeneous medium was widely used in geomechanics. A serious reason to reconsider the conventional views has been given by investigations of recent years, which indicate the necessity to take into account the block structure of mountain rocks in mathematical models intended for geomechanics and seismology [1, 2]. In this case, the mountain massif is considered as a system of nested blocks of various scale levels connected by interlayers consisting of weaker cracked rocks. The presence of such liable interlayers leads to the case when the block massif is deformed mainly due to their deformation both statically and dynamically.

In the simplest case, the dynamics of the block medium is investigated approximately considering that such blocks are incompressible and all deformations and displacements occur due to compressibility of interlayers [3–5]. The computing model can be a lattice of masses connected with each other by springs.

In terms of this model, we consider the anti-plane deformation of a two-dimensional square lattice consisting of masses $M$ connected by springs with length $L$ having identical rigidities $k$ in both directions. The unsteady propagation of perturbations under the effect of a step load is investigated. The problem in this statement with the sinusoidal load having the resonant frequency is solved analytically in [6].

The equations of motion of masses have the following form:

$$M\ddot{u}_{m,n} = k(u_{m+1,n} + u_{m-1,n} + u_{m,n+1} + u_{m,n-1} - 4u_{m,n}) + Q(t), \qquad (1)$$

where $u_{m,n}$ is the displacement of masses in the direction orthogonal to the lattice plane, $m$ and $n$ are the numbers of masses in the direction of axes $x$ and $y$, $Q(t) = Q_0 H(t)\delta_{m0}\delta_{n0}$ is a step load applied in the point with coordinates (0, 0), $H(t)$ is the Heaviside function, and $\delta_{n0}$ is the Kronecker delta.

Applying the discrete Fourier transformation over variables $m$ and $n$ and the Laplace transformation over time, we derive the solution in images**:**

$$u^{LF_mF_n} = \frac{Q^{LF_mF_n}}{Mp^2 + 2k[2-\cos(q_xL)-\cos(q_yL)]}, \quad Q^{LF_mF_n} = \frac{Q_0}{p}. \tag{2}$$

Here, $L$ denotes the Laplace transformation over time with parameter $p$; $F_m$ and $F_n$ are the discrete Fourier transformations over $m$ and $n$ with parameters $q_x$ and $q_y$, respectively.

Inverting derived expression (2) over variable $q_y$, we obtain

$$u_n^{LF_m} = \frac{Q_0\left(B-\sqrt{B^2-1}\right)^{|n|}}{2kp\sqrt{B^2-1}}, \quad B = \frac{Mp^2}{2k} + 2 - \cos(q_xL).$$

Let us investigate the unsteady behavior of long-wave perturbations at large time intervals. We assume that $p \to 0$. This corresponds to $t \to \infty$ in the space of the originals. The asymptotics of images of displacements at $n=0$ will be as follows:

$$u_0^{LF_m} \approx \frac{Q_0}{2kp\sqrt{B^2-1}}, \quad B^2-1 \approx 4\left(\frac{p^2}{4\omega_0^2} + \sin^2\frac{q_xL}{2}\right)\left(1+\sin^2\frac{q_xL}{2}\right), \quad \omega_0 = \sqrt{\frac{k}{M}}. \tag{3}$$

Here and below, formula $v(z) \approx w(z)$ as $z \to z_0$ means that $\lim_{z \to z_0}[v(z)-w(z)] = 0$.

Using formulas of inversion of the Fourier and Laplace transformations [7–9], we derive the asymptotics of long-wave perturbations as $t \to \infty$:

$$u_{m,0} \approx \frac{Q_0}{2k\pi}\begin{cases} H(z-1)\ln\left(z+\sqrt{z^2-1}\right), & z = \omega_0 t/|m|, \quad m \neq 0, \\ \ln(4\sqrt{2}\omega_0 t) + \gamma, & m = 0, \end{cases} \tag{4}$$

$$\dot{u}_{m,0} \approx \frac{Q_0}{2\sqrt{kM}} J_m^2(\omega_0 t), \tag{5}$$

$$\ddot{u}_{m,0} \approx \frac{Q_0}{2M} J_m(\omega_0 t)\left[J_{m-1}(\omega_0 t) - J_{m+1}(\omega_0 t)\right]. \tag{6}$$

Here, $J_m$ is a Bessel function of the first kind of order $m$ [8] and $\gamma = 0.577...$ is the Euler constant.

Using the approximate representation of the Bessel functions through the Airy function, which is valid at $m \gg 1$ [8, 9], let us rewrite solutions (5) and (6) in another form:

$$\dot{u}_{m,0} \approx \frac{Q_0}{2^{1/3}\sqrt{kM}}\left[\frac{\text{Ai}(2^{1/3}\kappa)}{(\omega_0 t)^{1/3}}\right]^2, \quad \kappa = \frac{m-\omega_0 t}{(\omega_0 t)^{1/3}}, \quad t \to \infty. \tag{7}$$

$$\ddot{u}_{m,0} \approx -\frac{2Q_0\text{Ai}(2^{1/3}\kappa)\text{Ai}'(2^{1/3}\kappa)}{t\sqrt{kM}}, \quad t \to \infty. \tag{8}$$

Here, Ai is the Airy function [9], and the prime means the derivative with respect to the argument.

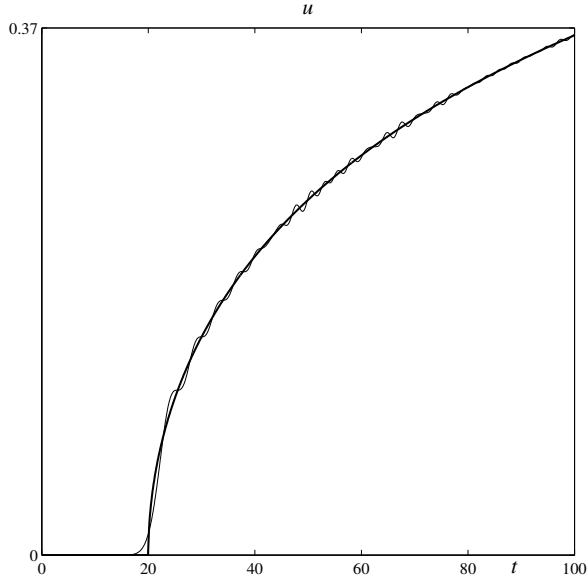 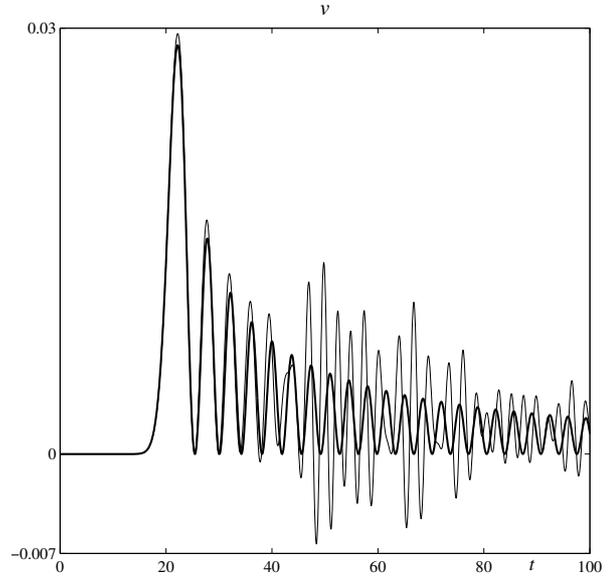

**Fig. 1.** Oscillograms of the displacements  **Fig. 2**. Oscillograms of the velocity of displacements

The advantage of solutions (7) and (8) compared with (5) and (6) is in the fact that these solutions explicitly describe the degree of decay of long-wave perturbations near the quasi-front $mL = c_* t$, where $c_* = L\omega_0$ is the velocity of infinitely long waves in the lattice. The analysis of solutions (7) and (8) shows that the maximal amplitudes of the velocity of displacement drop as $t^{-2/3}$ and those of accelerations drop as $t^{-1}$ as $t \to \infty$. The quasi-front zone, where perturbations vary from zero to maximum, is extended as $t^{1/3}$ as $t \to \infty$.

Equation (1) was also solved by the finite difference method using the explicit scheme

$$u_{m,n}^{r+1} - 2u_{m,n}^r + u_{m,n}^{r-1} = \tau^2 \left( u_{m+1,n}^r + u_{m-1,n}^r + u_{m,n+1}^r + u_{m,n-1}^r - 4u_{m,n}^r + Q_0 \delta_{m0} \delta_{n0} \right), \ r \geq 0. \qquad (9)$$

Here, $t = r\tau$, where $\tau$ is the time step of the difference mesh, and $r$ is the layer number over time. Figs. 1–3 show oscillograms of displacements $u = u_{m,0}$ (Fig. 1), velocities of displacements $v = \dot{u}_{m,0}$ (Fig. 2), and accelerations $w = \ddot{u}_{m,0}$ (Fig. 3) under the step load. The parameters of the problem are as follows: $m = 20$, $k = 1$, $L = 1$, $M = 1$ and $Q_0 = 1$. Thin curves correspond to the finite difference solution by scheme (9) ($\tau = 0.07$), and thick curves correspond to analytical solutions (4), (7), and (8). The comparison of numerical and asymptotic solutions (4), (7), and (8) shows that the long-wave asymptotics accurately describes the solution near the quasi-front $m = c_* t$. Starting from moment of time $t = m\sqrt{2}/c \sim 44$, where $c = 2/\pi$ is the phase velocity of short-wave perturbations, perturbations appear in oscillograms of numerical solutions of $v$ and $w$, which correspond to oscillations of short waves ($q_x = q_y = \pi$); the latter are absent in asymptotic solutions (7) and (8). We also compared the finite-difference solutions for velocities of displacements $v$ and accelerations $w$ with analytical solutions (5) and (6). As should be expected, their coincidence is even better than with solutions (7) and (8).

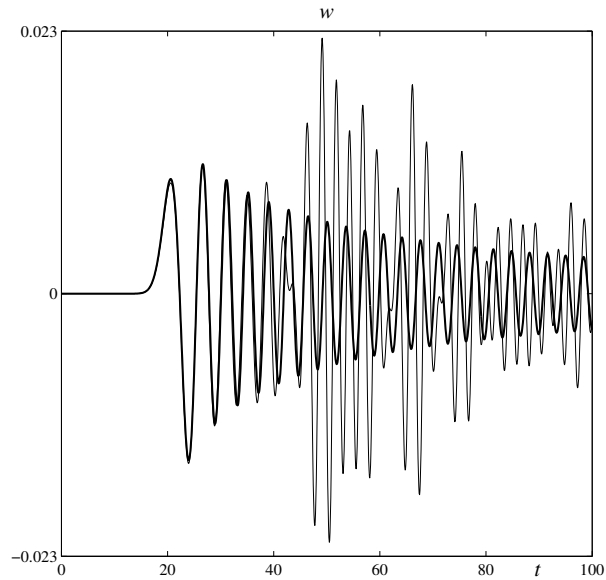

**Fig. 3.** Oscillograms of the accelerations

Thus, under the step impact on the square lattice of spring-connected masses, the following effects are observed in the lattice as $t \to \infty$:

(i) the amplitudes of displacements grow proportionally to $\ln\left(z + \sqrt{z^2 - 1}\right)$, where $z = \omega_0 t/|m|$, $m \neq 0$;

(ii) the maximal amplitudes of velocities of displacements of long-wave disturbances drop with time as $t^{-2/3}$, and those of accelerations as $t^{-1}$;

(iii) the quasi-front zone is extended as $t^{1/3}$.

This comparison of the numerical and analytical solutions shows that the asymptotic solution derived as $t \to \infty$ accurately describes long-wave perturbations. The agreement of these solutions appears at a finite impact time or at a finite distance from the place of applying the load.